\documentclass[11pt]{amsart}
\usepackage[utf8]{inputenc}
\usepackage[T1]{fontenc}
\usepackage{amsmath,amssymb,amsthm}
\usepackage[margin=2.5cm]{geometry}
\usepackage{tikz-cd}
\usepackage{hyperref}
\hypersetup{
  pdftitle={Separately continuous maps into regular spaces and continuity of operations in separable semitopological groups},
  pdfauthor={Evgenii Reznichenko},
  pdfsubject={General Topology (MSC2020: 54H11, 22A20, 54C08, 54E52, 91A44)},
  pdfkeywords={semitopological group, feebly compact space, separately continuous map, strong quasicontinuity, Banach-Mazur game, Baire space}
}

\newtheorem{theorem}{Theorem}[section]
\newtheorem{proposition}[theorem]{Proposition}
\newtheorem{lemma}[theorem]{Lemma}
\newtheorem{corollary}[theorem]{Corollary}
\theoremstyle{definition}
\newtheorem{definition}[theorem]{Definition}
\newtheorem{question}{Question}
\theoremstyle{remark}
\newtheorem{remark}[theorem]{Remark}

\newcommand{\R}{\mathbb{R}}
\newcommand{\om}{\omega}
\newcommand{\bt}{\beta}
\newcommand{\cl}[1]{\overline{#1}}
\newcommand{\Int}{\operatorname{Int}}
\newcommand{\Ls}{\operatorname{Ls}}
\newcommand{\term}[1]{\emph{#1}}

\title[Separately continuous maps and separable semitopological groups]{Separately continuous maps into regular spaces and continuity of operations in separable semitopological groups}
\author{Evgenii Reznichenko}
\address{Department of General Topology and Geometry, Faculty of Mechanics and Mathematics\\
Lomonosov Moscow State University, Leninskie Gory 1, Moscow 119991, Russia}
\address{Shenzhen MSU-BIT University, 1 International University Park Road, Dayun New Town, Longgang District\\ Shenzhen, Guangdong 518172, China}
\email{erezn@inbox.ru}
\subjclass[2020]{54H11, 22A20, 54C08, 54E52, 91A44}
\keywords{Semitopological group, feebly compact space, separately continuous map, strong quasicontinuity,
Banach--Mazur game, Baire space}
\date{}

\begin{document}
\begin{abstract}
We prove that every regular separable feebly compact semitopological group is a topological group. This answers a question of M. Tkachenko (2014) and extends the separable case of the author's 1994 theorem on Tychonoff pseudocompact groups to regular spaces. The proof rests on a new result: if $X$ is a separable Baire space, $Y$ is feebly compact and $Z$ is regular, then every separately continuous map $X\times Y\to Z$ is strongly quasicontinuous. Previously this was known only when $X$, $Y$ and $Z$ are Tychonoff. Using topological games of Banach--Mazur type, we extend the group theorem to regular separable Baire groups from a wider class of spaces, which contains all $p$-spaces, all feebly compact spaces and all products of a $p$-space with a feebly compact space. We also obtain local versions: the theorem holds for locally separable, locally feebly compact groups.
\end{abstract}

\maketitle

\section{Introduction}

All spaces are topological spaces; no separation axioms are assumed unless they are stated explicitly. Regular spaces are assumed to be
$T_1$-spaces, and Tychonoff spaces are completely regular $T_1$-spaces. A map
$f\colon X\times Y\to Z$ is \term{separately continuous} if the maps $x\mapsto f(x,y_*)$ and $y\mapsto f(x_*,y)$ are continuous for all
$x_*\in X$ and $y_*\in Y$. A group $G$ with a topology is \term{semitopological} if the multiplication $(x,y)\mapsto xy$ is separately
continuous, \term{paratopological} if the multiplication is continuous, and \term{topological} if the multiplication and the inversion
$x\mapsto x^{-1}$ are continuous.

\subsection*{The problem}
The starting point of this paper is the following theorem.

\begin{theorem}[{\cite[Corollary~2.7]{R1994}}]\label{th-R1994}
Let $G$ be a Tychonoff pseudocompact semitopological group. If $G$ is countably compact, or has countable tightness, or is separable,
or is a $k$-space, then $G$ is a topological group.
\end{theorem}

Recall that a Tychonoff space is \term{pseudocompact} if every continuous real-valued function on it is bounded \cite{Hewitt1948}.
Beyond the Tychonoff setting the natural analogue is the following: a space is \term{feebly compact} if every locally finite family of
nonempty open subsets is finite. A Tychonoff space is feebly compact if and only if it is pseudocompact \cite[Theorem~3.10.22]{Engelking1989}.

Every $T_0$ topological group is Tychonoff \cite{Pontryagin1939}, and, by a remarkable theorem of Banakh and Ravsky \cite{BR2017},
so is every regular paratopological group. This does not extend to semitopological groups:
there exists an infinite regular feebly compact Boolean semitopological group (that is, of period~$2$) on which every continuous real-valued function is constant \cite{RT2026};
in particular, it is not Tychonoff. For this reason the proof in \cite{R1994}, which extends the multiplication from $G\times G$ to
$\bt G\times\bt G$, where $\bt G$ is the Stone--\v{C}ech compactification of $G$, does not apply to regular spaces, and the difference between regular and
Tychonoff spaces is essential here.

In the countably compact case the regular version is nevertheless known:
as observed in \cite[p.~838]{T2014}, it follows from \cite{KKM2001}, since regular countably compact spaces are strongly Baire and
every strongly Baire semitopological group is topological. Tkachenko asked whether ``Tychonoff pseudocompact'' in
Theorem~\ref{th-R1994} can be weakened to ``regular feebly compact''.

\begin{question}[{\cite[Problem~3.19]{T2014}}]\label{q-T}
Let $G$ be a regular feebly compact semitopological group. Is $G$ a topological group provided that
(a)~$G$ is separable, (b)~$G$ has countable tightness, or (c)~$G$ is a $k$-space?
\end{question}

We give an affirmative answer to part~(a).

\begin{theorem}\label{th-main-group}
Every regular separable feebly compact semitopological group is a topological group.
\end{theorem}

Since $T_1$ topological groups are Tychonoff, every group as in Theorem~\ref{th-main-group} is Tychonoff and hence pseudocompact.
Thus separability rules out non-Tychonoff groups such as the one constructed in \cite{RT2026}.

\subsection*{The strategy}
A map $g\colon T\to Z$ is \term{quasicontinuous at} $t\in T$ if for every neighborhood $W$ of $t$ and every neighborhood $O$ of $g(t)$
there is a nonempty open set $V\subset W$ with $g(V)\subset O$. For a map $f\colon X\times Y\to Z$ this means that for all neighborhoods
$U$ of $x_*$, $V$ of $y_*$ and $O$ of $f(x_*,y_*)$ there are nonempty open sets $U'\subset U$ and $V'\subset V$ with $f(U'\times V')\subset O$.
We say that $f$ is \term{strongly quasicontinuous at $(x_*,y_*)$} if, in addition, $V'$ can always be chosen to be a neighborhood of $y_*$
\cite{Bouziad1996,KKM2001,R2024FA}. For the multiplication of a semitopological group quasicontinuity and strong quasicontinuity
are equivalent, and both are equivalent to feeble continuity (Remark~\ref{rem-qc-types}).

A set $P\subset X\times X$ is a \term{semi-neighborhood of the diagonal} if $x\in\Int\{y:(x,y)\in P\}$ for every $x\in X$.
A space $X$ is \term{$\Delta$-nonmeager} \cite{R2025} if for every semi-neighborhood $P$ of the diagonal there is a nonempty open set
$W\subset X$ such that $P\cap(W\times W)$ is dense in $W\times W$.\footnote{These spaces were called $\Delta$-Baire in
\cite{R2022ci,Moors2017}. For homogeneous spaces, in particular for groups, the two names refer to the same class
\cite[Theorem~3.1(5)]{R2025}.}

\begin{theorem}[{\cite[Theorem~1]{Moors2017}; see also \cite[Theorem~7.8]{R2024}}]\label{th-delta-stg}
Let $G$ be a regular $\Delta$-nonmeager semitopological group. If the multiplication on $G$ is quasicontinuous, then $G$ is
a topological group.
\end{theorem}

The assumption that $G$ is $\Delta$-nonmeager cannot be dropped: the Sorgenfrey line is a regular Baire paratopological group,
but not a topological group, hence it is not $\Delta$-nonmeager \cite[Example~9.4]{R2025}.
Theorem~\ref{th-delta-stg} suggests a scheme for proving that a semitopological group $G$ is topological:
(1)~prove that $G$ is $\Delta$-nonmeager; (2)~prove that the multiplication on $G$ is quasicontinuous. As noted in \cite{R2024}, nearly all
known results on the continuity of operations in semitopological groups, apart from the paratopological case, fit this scheme;
a rare exception is Theorem~\ref{th-R1994}.

For a regular feebly compact group, step~(1) is known \cite{R2022ci,R2025,Moors2017,R2024}. Step~(2), however, may fail:
Korovin \cite{Korovin1992} constructed a pseudocompact Abelian semitopological group of period~$2$ which is not a topological group
(see also \cite{R1994}); this group is regular and feebly compact, so by Theorem~\ref{th-delta-stg} its multiplication is not quasicontinuous.
Groups of this kind are very far from being separable: in the pseudocompact quasitopological groups obtained by Korovin's method every
countable subset is closed, discrete and $C^*$-embedded \cite{RT2024}. Thus the whole point is to show that separability forces
the multiplication to be quasicontinuous.

\subsection*{Separately continuous maps}
Quasicontinuity of the multiplication is a special case of the following classical question, which goes back to Baire; see the monograph \cite{CM2024}.

\begin{question}\label{q-sc}
Under which conditions on $X$, $Y$ and $Z$ is every separately continuous map $f\colon X\times Y\to Z$ quasicontinuous?
\end{question}

The case $Z=\R$ is the most thoroughly studied, and it covers all Tychonoff spaces $Z$. Regular spaces $Z$ were considered in
\cite{KKM2001} and \cite[Section~2.2]{CM2024}. We obtain a new sufficient condition, in which the second factor must satisfy
the following feeble analogue of the $q$-point property.

\begin{definition}\label{def-fq}
Let $X$ be a space. A sequence $(U_n)_{n\in\om}$ of nonempty open subsets of $X$ is \term{feebly compact} if for every choice of nonempty open
sets $V_n\subset U_n$, $n\in\om$, the sequence $(V_n)$ has an accumulation point in $X$, that is, a point every neighborhood of which
meets $V_n$ for infinitely many $n$. A point $x\in X$ is a \term{feebly $q$-point} if it has a feebly compact sequence of open
neighborhoods; $X$ is a \term{feebly $q$-space} if all its points are feebly $q$-points.
\end{definition}

A space is feebly compact if and only if every sequence of its nonempty open subsets is feebly compact; in particular, every feebly
compact space is a feebly $q$-space (Proposition~\ref{prop-fc}). The discrete space $\mathbb N$ is a feebly $q$-space which is not feebly compact.

\begin{theorem}\label{th-main-scf}
Let $X$ be a separable Baire space, $Y$ a feebly $q$-space and $Z$ a regular space. Then every separately continuous map
$f\colon X\times Y\to Z$ is strongly quasicontinuous.
\end{theorem}

For Tychonoff spaces the conclusion of Theorem~\ref{th-main-scf} can be derived from known results (see Remark~\ref{rem-tychonoff}), but all such
arguments use real-valued functions on $Z$ and do not extend to regular spaces, which is essential for groups. For regular ranges the closest known results
are the following. By the theorem of Martin \cite{Martin1961}, extended to regular ranges by
Neubrunn \cite{Neubrunn1976}, if $X$ is a Baire space, $Y$ is second countable and $Z$ is regular, then every separately continuous map
$X\times Y\to Z$ is quasicontinuous. Here the Baire property and the countability assumption concern different factors, whereas in
Theorem~\ref{th-main-scf} both are imposed on $X$, and the countability assumption is separability, which is much weaker than second countability.
Structurally, Theorem~\ref{th-main-scf} is closer to a theorem of Bouziad \cite[Theorem~1]{Bouziad1996}: a separately continuous map
$f\colon X\times Y\to Z$ is strongly quasicontinuous at every point of $X\times\{b\}$ provided that $X$ is a pointwise countably complete
Baire space (for example, a Baire $p$-space), $Z$ is completely regular and $b$ is a $q$-point of $Y$. From it Bouziad deduced that every
semitopological group which is a Baire $p$-space has continuous multiplication. In Theorem~\ref{th-main-scf} the completeness-type
assumption on $X$ is replaced by separability, the $q$-point property of $Y$ is weakened to its feeble analogue, and $Z$ may be an arbitrary regular space.

Our proof uses only the Banach--Mazur game on $X$: the feebly $q$-point property provides accumulation points of sequences of open sets in $Y$,
separability of $X$ reduces the control of $f$ to countably many conditions, and the Baire property of $X$, through Oxtoby's theorem, produces
the point of $X$ which leads to a contradiction. Applied to the multiplication $G\times G\to G$, Theorem~\ref{th-main-scf} yields step~(2)
of the scheme above, and Theorem~\ref{th-main-group} follows.

\subsection*{Wider classes of groups}
The proof of Theorem~\ref{th-main-group} uses feeble compactness only through two consequences: the group is
$\Delta$-nonmeager (step~(1)), and it is a Baire space in which every point is a feebly $q$-point (step~(2)). Both hold for a much wider class of spaces,
the \term{$\Gamma^{BM}_o$-nonmeager} spaces, defined by a modification of the Banach--Mazur game \cite{R2023} (Definition~\ref{def-gamma});
for regular spaces, the closely related $\Gamma^{BM}_o$-Baire spaces are exactly those considered in \cite{ACK2010} and \cite{R2022ci}
(Remark~\ref{rem-ACK}). Arhangel'skii, Choban and Kenderov proved that every regular paratopological group which is $\Gamma^{BM}_o$-Baire is a
topological group \cite[Corollary~3.4]{ACK2010}. For semitopological groups we prove the following.

\begin{theorem}\label{th-main-gamma}
Every regular separable $\Gamma^{BM}_o$-nonmeager semitopological group is a topological group.
\end{theorem}

Regular feebly compact spaces are $\Gamma^{BM}_o$-nonmeager (Lemma~\ref{l-baire}, Propositions~\ref{prop-classes}(1) and~\ref{prop-gamma-space}, and \cite[Propositions~7 and~11]{R2023}), so Theorem~\ref{th-main-group} is a special case of Theorem~\ref{th-main-gamma}.
The separability assumption cannot be dropped: Korovin's group is regular and feebly compact, hence $\Gamma^{BM}_o$-nonmeager, but it is not a topological group.

A natural source of $\Gamma^{BM}_o$-nonmeager spaces is the following feeble analogue of Arhangel'skii's $p$-spaces.
It contains all $p$-spaces, more generally all pointwise countably complete spaces of Bouziad, all feebly compact spaces, and the products of
a $p$-space with a feebly compact space (Section~\ref{sec-fp-vs-p}).

\begin{definition}\label{def-fp}
A space $X$ is a \term{feebly $p$-space} if there is a sequence $(\gamma_n)_{n\in\om}$ of open covers of $X$ with the following property:
if $x\in X$, $x\in U_n\in\gamma_n$ for $n\in\om$, and $V_n\subset U_n$ are nonempty open sets with $V_{n+1}\subset V_n$, then the sequence
$(V_n)$ is feebly compact.
\end{definition}

Every regular Baire feebly $p$-space is $\Gamma^{BM}_o$-nonmeager (Proposition~\ref{prop-gamma-space} and \cite[Propositions~7 and~11]{R2023}),
which gives the following two results.

\begin{proposition}\label{prop-fp-delta}
Every regular Baire feebly $p$-space is $\Delta$-nonmeager.
\end{proposition}

\begin{theorem}\label{th-main-pspace}
Let $G$ be a regular separable Baire semitopological group. If $G$ is a feebly $p$-space, then $G$ is a topological group.
\end{theorem}

For groups which are regular Baire $p$-spaces the conclusion of Theorem~\ref{th-main-pspace} is known even without separability, by \cite{Bouziad1996} and \cite{ACK2010};
the new cases are, for example, feebly compact groups and groups homeomorphic to products of a $p$-space with a feebly compact space.
By Proposition~\ref{prop-fp-delta}, the Sorgenfrey line, being a regular Baire space which is not $\Delta$-nonmeager, is not a feebly $p$-space.

\subsection*{Local versions}
All the main results remain valid, with the same proofs, for spaces satisfying the assumptions locally (Remark~\ref{rem-local}): in
Theorems~\ref{th-main-group}, \ref{th-main-scf}, \ref{th-main-gamma} and~\ref{th-main-pspace} separability can be replaced by local separability, and in
Proposition~\ref{prop-fp-delta} and Theorem~\ref{th-main-pspace} feebly $p$-spaces can be replaced by locally feebly $p$-spaces.
In particular, every regular locally separable locally feebly compact semitopological group is a topological group (Corollary~\ref{cor-local}).

\subsection*{Organization}
Section~\ref{sec-prelim} contains preliminaries. Theorem~\ref{th-main-scf} is proved in Section~\ref{sec-qc}, and
Theorem~\ref{th-main-group} in Section~\ref{sec-group}; these sections suffice for the answer to Question~\ref{q-T}. Section~\ref{sec-fp}
is devoted to feebly $p$-spaces and $\Gamma^{BM}_o$-nonmeager spaces and contains the proofs of the remaining results and the local versions.
Open problems are collected in Section~\ref{sec-problems}.

\section{Preliminaries}\label{sec-prelim}

\subsection{Accumulation points of sequences of sets}\label{sec-acc}
Let $(A_n)_{n\in\om}$ be a sequence of subsets of a space $X$. A point $p\in X$ is an \term{accumulation point} of $(A_n)$, and $(A_n)$
\term{accumulates} to $p$, if every neighborhood of $p$ meets $A_n$ for infinitely many $n$. The set of accumulation points of $(A_n)$ is
\[
\Ls_nA_n=\bigcap_{n\in\om}\cl{\bigcup_{k\ge n}A_k}.
\]
If $A_n=\{x_n\}$ for all $n$, then $\Ls_nA_n$ is the set of cluster points of the sequence $(x_n)$. We shall often use the following
obvious fact:
\begin{equation}\label{eq-ls}
\text{if } A_k\subset B \text{ for all } k\ge m, \text{ then } \Ls_nA_n\subset\cl{B}.
\end{equation}

\subsection{Feebly compact spaces}
Lemmas~\ref{l-fc} and~\ref{l-baire} are well known for pseudocompact spaces (see \cite[Theorem~1.1.3(4) and Theorem~1.2.2]{PTS2018}; the
latter is due to Colmez), and Lemma~\ref{l-fc} is stated for feebly compact spaces in \cite[p.~64]{Tk2018}. We include the short proofs for
completeness.

\begin{lemma}\label{l-fc}
Let $X$ be a feebly compact space and let $O_0\supset O_1\supset\dots$ be nonempty open subsets of $X$.
Then $\bigcap_{i}\cl{O_i}\ne\emptyset$.
\end{lemma}

\begin{proof}
Suppose that $\bigcap_i\cl{O_i}=\emptyset$ and put ${E}_i=O_i\setminus\cl{O_{i+1}}$. The sets ${E}_i$ are open and pairwise disjoint: if
$i<j$, then ${E}_j\subset O_j\subset O_{i+1}$, while ${E}_i\cap O_{i+1}=\emptyset$. If ${E}_i=\emptyset$ for all $i\ge i_0$, then
$O_i\subset\cl{O_{i+1}}$, hence $\cl{O_i}=\cl{O_{i+1}}$ for $i\ge i_0$ and $\bigcap_i\cl{O_i}=\cl{O_{i_0}}\ne\emptyset$, a contradiction.
Thus infinitely many sets ${E}_i$ are nonempty. The family $({E}_i)$ is locally finite: given $x\in X$, choose $j$ with $x\notin\cl{O_j}$;
then the neighborhood $X\setminus\cl{O_j}$ of $x$ misses ${E}_i\subset O_j$ for all $i\ge j$. This contradicts the feeble compactness
of $X$.
\end{proof}

\begin{lemma}\label{l-baire}
Every regular feebly compact space is a Baire space.
\end{lemma}

\begin{proof}
Let ${E}_1,{E}_2,\dots$ be open dense subsets of a regular feebly compact space $X$ and let $W_0$ be a nonempty open set. By regularity,
there are nonempty open sets $W_i$ with $\cl{W_i}\subset W_{i-1}\cap {E}_i$, $i\ge1$. By Lemma~\ref{l-fc}, $\bigcap_i\cl{W_i}\ne\emptyset$,
and this set is contained in $W_0\cap\bigcap_i{E}_i$.
\end{proof}

\begin{proposition}\label{prop-fc}
A space $X$ is feebly compact if and only if every sequence of nonempty open subsets of $X$ is feebly compact. In particular, every
point of a feebly compact space is a feebly $q$-point.
\end{proposition}

\begin{proof}
Let $X$ be feebly compact and let $M_n$, $n\in\om$, be nonempty open subsets of $X$. The sets $T_n=\bigcup_{k\ge n}M_k$ are nonempty, open
and decreasing, so by Lemma~\ref{l-fc} there is a point $p\in\bigcap_n\cl{T_n}=\Ls_nM_n$. Hence every sequence of nonempty open sets is
feebly compact; in particular, the sequence $(X,X,\dots)$ of neighborhoods of any point is feebly compact. Conversely, if $X$ is not
feebly compact, choose pairwise distinct members $N_k$, $k\in\om$, of an infinite locally finite family of nonempty open sets. The
sequence $(N_k)$ has no accumulation point, so the sequence $(X,X,\dots)$ is not feebly compact.
\end{proof}

\subsection{The Banach--Mazur game}\label{sec-bm}
Let $X$ be a space. In the \term{Banach--Mazur game} on $X$ two players, $\alpha$ and $\beta$, choose nonempty open subsets of $X$.
We put $U_0=X$, and the game consists of steps $n=1,2,\dots$:

\smallskip
\emph{Step $n$.} Player $\beta$ chooses a nonempty open set $V_n\subset U_{n-1}$; then player $\alpha$ chooses a nonempty open set
$U_n\subset V_n$.

\smallskip
\noindent A \term{play} is the resulting sequence of moves. Since $V_{n+1}\subset U_n\subset V_n$, the sets $V_n$ decrease and
$\bigcap_nV_n=\bigcap_nU_n$. A \term{strategy} of a player assigns to every finite sequence of previous moves an admissible move of this
player; a play \term{follows} a strategy if all moves of this player in the play are given by the strategy.

\begin{theorem}[Oxtoby \cite{Oxtoby1957}]\label{th-oxtoby}
A space $X$ is a Baire space if and only if $\beta$ has no strategy in the Banach--Mazur game on $X$ such that $\bigcap_nV_n=\emptyset$
in every play which follows it.
\end{theorem}

We shall use only the implication from left to right.

\section{Strong quasicontinuity of separately continuous maps}\label{sec-qc}

Theorem~\ref{th-main-scf} is an immediate consequence of the following pointwise statement, since every point of a feebly $q$-space is a
feebly $q$-point.

\begin{theorem}\label{th-main-scf+}
Let $X$ be a separable Baire space, $Y$ a space and $Z$ a regular space, and let $f\colon X\times Y\to Z$ be a separately continuous map.
If $y_*\in Y$ is a feebly $q$-point, then $f$ is strongly quasicontinuous at $(x_*,y_*)$ for every $x_*\in X$.
\end{theorem}

\begin{proof}
Let $z_*=f(x_*,y_*)$, and suppose that $f$ is not strongly quasicontinuous at $(x_*,y_*)$. Then there are open sets $U_*'\ni x_*$,
$V_*\ni y_*$ and $O\ni z_*$ such that
\begin{equation}\label{eq-fail}
\begin{gathered}
f(\tilde U\times\tilde V)\not\subset O\quad\text{for every nonempty open }\tilde U\subset U_*'\\
\text{and every open }\tilde V\text{ with }y_*\in\tilde V\subset V_*.
\end{gathered}
\end{equation}

\emph{Step 1: the sets $A$, $B$, $U_*$ and $D$.} Since $Z$ is regular, there are open sets $A$ and $A'$ with $z_*\in A$,
$\cl{A}\subset A'$ and $\cl{A'}\subset O$. Put $B=Z\setminus\cl{A'}$. Since $f(\cdot,y_*)$ is continuous, there is an open neighborhood
$U_*\subset U_*'$ of $x_*$ with $f(U_*\times\{y_*\})\subset A$. Since $X$ is separable, so is its open subspace $U_*$; let
$D=\{x_n:n\in\om\}$ be a countable dense subset of $U_*$ (repetitions are allowed).

\emph{Step 2: the neighborhoods $W_n$.} Let $(Q_n)_{n\in\om}$ be a feebly compact sequence of open neighborhoods of $y_*$. Since
$f(x_n,y_*)\in A$ and the maps $f(x_n,\cdot)$ are continuous, we can choose inductively open neighborhoods $W_n$ of $y_*$ such that
\[
W_0\subset V_*,\qquad W_{n+1}\subset W_n,\qquad W_n\subset Q_n,\qquad f(\{x_n\}\times W_n)\subset A\qquad(n\in\om).
\]
Put ${C}=\bigcap_n\cl{W_n}$. We claim:
\begin{enumerate}
\item[(a)] $f(U_*\times{C})\subset\cl{A}$;
\item[(b)] if $S_n\subset W_n$ are nonempty open sets, $n\in\om$, then $(S_n)$ accumulates to a point of ${C}$;
\item[(c)] $\cl{A}\cap\cl{B}=\emptyset$;
\item[(d)] $f(\tilde U\times\tilde V)\cap B\ne\emptyset$ for every nonempty open $\tilde U\subset U_*$ and every open $\tilde V$ with
$y_*\in\tilde V\subset V_*$.
\end{enumerate}

(a) Let ${c}\in{C}$ and $n\in\om$. Since ${c}\in\cl{W_n}$, $f(\{x_n\}\times W_n)\subset A$ and $f(x_n,\cdot)$ is continuous, we have
$f(x_n,{c})\in\cl{A}$. Since $f(\cdot,{c})$ is continuous and $D$ is dense in $U_*$, it follows that $f(U_*\times\{{c}\})\subset\cl{A}$.

(b) The sets $S_n$ are nonempty open subsets of $Q_n$, and $(Q_n)$ is feebly compact, so $(S_n)$ accumulates to a point $p\in Y$. For every
$n$ we have $S_k\subset W_n$ for $k\ge n$, so $p\in\cl{W_n}$ by \eqref{eq-ls}. Hence $p\in{C}$.

(c) The set $Z\setminus A'$ is closed and contains $B$, so $\cl{B}\subset Z\setminus A'$, while $\cl{A}\subset A'$.

(d) Since $\tilde U\subset U_*'$, we have $f(\tilde U\times\tilde V)\not\subset O$ by \eqref{eq-fail}, and $Z\setminus O\subset B$ because
$\cl{A'}\subset O$.

\emph{Step 3: a strategy for $\beta$.} We define a strategy for $\beta$ in the Banach--Mazur game on $X$; all its moves lie in $U_*$. (The moves $V_n$ and $U_n$ of this game are subsets of $X$, whereas $V_*$ and the sets $W_n$ are subsets of $Y$.)
Let $n\ge0$, and suppose that the sets $V_1\supset U_1\supset\dots\supset V_n\supset U_n$ have been chosen; for $n=0$ we use $U_*$ in place
of $U_0$. By (d), applied to $\tilde U=U_n$ (respectively, $\tilde U=U_*$ if $n=0$) and $\tilde V=W_n$, there are $x'\in\tilde U$ and
$y_n\in W_n$ with $f(x',y_n)\in B$. Since $f(\cdot,y_n)$ is continuous and $B$ is open, there is a nonempty open set $V_{n+1}\subset\tilde U$
with
\[
f(V_{n+1}\times\{y_n\})\subset B;
\]
player $\beta$ chooses this set. (For each $n$ and each $U_n$, the point $y_n$ and the set $V_{n+1}$ are fixed once and for all.) Since
$X$ is a Baire space, Theorem~\ref{th-oxtoby} gives a play which follows this strategy and in which $\bigcap_{n\ge1}U_n\ne\emptyset$.
Fix $x\in\bigcap_{n\ge1}U_n$. Then $x\in V_{n+1}$ for all $n\in\om$; hence $x\in U_*$ and
\[
f(x,y_n)\in B\qquad\text{for all }n\in\om.
\]

\emph{Step 4: the contradiction.} The map $h=f(x,\cdot)$ is continuous and $B$ is open, so there are open sets $S_n$ with
$y_n\in S_n\subset W_n$ and $h(S_n)\subset B$. By (b), $(S_n)$ accumulates to a point $y\in{C}$. By the continuity of $h$, every
neighborhood of $h(y)$ meets $h(S_n)\subset B$ for infinitely many $n$, so $f(x,y)=h(y)\in\cl{B}$. On the other hand, $(x,y)\in U_*\times{C}$,
so $f(x,y)\in\cl{A}$ by (a). This contradicts (c).
\end{proof}

\begin{remark}
The feebly $q$-point property of $y_*$ is used only in (b). Separability of $X$ provides the countable set $D$, which together with the
neighborhoods $W_n$ gives (a). The Baire property of $X$ is used, through Oxtoby's theorem, to find the point $x$. Property (d) is applied
only to neighborhoods of $y_*$; this is why the argument gives strong quasicontinuity rather than quasicontinuity.
\end{remark}

\begin{remark}[The Tychonoff case]\label{rem-tychonoff}
For Tychonoff spaces $X$, $Y$, and $Z$, with $Y$ pseudocompact, the conclusion of Theorem~\ref{th-main-scf} can also be derived from
known results in at least two ways. Since a map into a Tychonoff space is strongly quasicontinuous if and only if all its compositions
with continuous functions $Z\to[0,1]$ are, it suffices to consider $Z=\R$; let $f\colon X\times Y\to\R$ be separately continuous.
In both ways one shows that $f$ has the Namioka property, that is, there is a dense $G_\delta$-set $A\subset X$ such that $f$ is continuous
at every point of $A\times Y$; this easily implies strong quasicontinuity (given $(x_*,y_*)$, choose $a\in A$ near $x_*$ with $f(a,y_*)$
near $f(x_*,y_*)$ and use the joint continuity of $f$ at $(a,y_*)$; cf.\ \cite[Proposition~1]{R2024FA}).
\begin{enumerate}
\item[(1)] Every separable Baire space is $\sigma$-$\beta$-unfavorable, that is, $\beta$ has no winning strategy in the
Christensen--Saint-Raymond game \cite{Christensen1981,SR1983} (see \cite[Section~3]{R2024FA} for the form used here): player $\alpha$ may
choose the points of a fixed countable dense set, and then Oxtoby's theorem applies. If $X$ is $\sigma$-$\beta$-unfavorable and $Y$ is
pseudocompact, then every separately continuous function $f\colon X\times Y\to\R$ has the Namioka property \cite[Proposition~3]{R2024FA}.
\item[(2)] Let $\phi\colon Y\to C_p(X)$, $\phi(y)(x)=f(x,y)$. The map $\phi$ is continuous, so $\phi(Y)$ is pseudocompact. Since $X$ is
separable, the restriction to a countable dense subset of $X$ is a continuous injection of $C_p(X)$ into a metrizable space, and its
restriction to the pseudocompact Tychonoff space $\phi(Y)\subset C_p(X)$ is a condensation onto a space with a countable base. A
condensation of a pseudocompact space onto a space with a countable base is a homeomorphism \cite[Proposition~II.6.2(c)]{Arh1992} (cf.\ \cite{AP1984}); hence $\phi(Y)$ is
metrizable, and therefore compact. By \cite[Assertion~1.2]{R1994} (with the factors interchanged), $f$ extends to a separately continuous function on $X\times\bt Y$. Separable
Baire spaces are Namioka spaces \cite[Theorem~4.3.3]{CM2024}; since $\bt Y$ is compact, the extension, and hence $f$, has the Namioka property.
\end{enumerate}
Both arguments are based on real-valued functions on $Z$, and the second one also on the Stone--\v{C}ech compactification of $Y$. They are
not available when $Z$ is regular but not Tychonoff, which is essential for groups: there are regular feebly compact Boolean semitopological groups
on which every continuous real-valued function is constant \cite{RT2026}. The proof of Theorem~\ref{th-main-scf} above is a version of the
game-theoretic argument of \cite[Theorem~3]{R2024FA}, which goes back to Christensen and Saint-Raymond, adapted to regular ranges.

Other related results: Calbrix and Troallic \cite{CT1979} proved that if $Y$ is second countable and $Z$ is metrizable, then for every space $X$
and every separately continuous map $f\colon X\times Y\to Z$ the set of points of discontinuity of $f$ has a meager projection on $X$
(see also \cite[Theorem~2.4.3]{CM2024}). For Tychonoff pseudocompact spaces $X$ and $Y$, quasicontinuity of a separately continuous function
$X\times Y\to\R$ is equivalent to its extendability to a separately continuous function on $\bt X\times Y$ \cite[Theorem~1]{R2022}.
Quasicontinuity was introduced by Kempisty \cite{Kempisty1932}.
\end{remark}

\section{\texorpdfstring{Proof of Theorem~\ref{th-main-group}}{Proof of the main theorem}}\label{sec-group}

Let $G$ be a regular separable feebly compact semitopological group. We follow the scheme described in the introduction.

(1) By Lemma~\ref{l-baire}, $G$ is a Baire space, and regular feebly compact spaces are $\Delta$-nonmeager (see \cite[Theorem~4.1]{R2024}
and \cite{R2022ci,R2025,Moors2017}; this also follows from Proposition~\ref{prop-fp-delta}). Hence $G$ is $\Delta$-nonmeager.

(2) The multiplication $f\colon G\times G\to G$, $f(x,y)=xy$, is separately continuous. By Proposition~\ref{prop-fc}, $G$ is a feebly
$q$-space; since $G$ is also a separable Baire space and a regular space, Theorem~\ref{th-main-scf} with $X=Y=Z=G$ shows that $f$ is
strongly quasicontinuous, and in particular quasicontinuous.

By Theorem~\ref{th-delta-stg}, $G$ is a topological group. \qed

\begin{remark}[Kinds of quasicontinuity]\label{rem-qc-types}
A map $g\colon T\to Z$ is \term{feebly continuous} if $\Int g^{-1}(O)\ne\emptyset$ for every open set $O\subset Z$ with
$g^{-1}(O)\ne\emptyset$. Strong quasicontinuity implies quasicontinuity, which implies feeble continuity. For the multiplication of a
semitopological group the three properties are equivalent \cite[Proposition~3.1(5) and footnote~3]{R2024}; in \cite{R2024} this is
stated for quasicontinuity with respect to the first variable, and the case of the second variable follows by passing to the opposite
multiplication $(x,y)\mapsto yx$. Strong quasicontinuity is called quasicontinuity with respect to the second variable in \cite{CM2024,R2024}.
By translations, it suffices to verify strong quasicontinuity of the multiplication at $(e,e)$. In \cite{Moors2017} the multiplication is assumed to be feebly continuous,
which is equivalent to quasicontinuity by the above; \cite{R2024} contains several strengthenings of Theorem~\ref{th-delta-stg}.
For general groups with a topology the three properties differ \cite[Example~2.6 and Section~3]{CDP2010}. Classes of $\Delta$-nonmeager spaces are studied in
\cite{R2022ci,R2023,R2025,Moors2017}.
\end{remark}

\section{\texorpdfstring{Feebly $p$-spaces}{Feebly p-spaces}}\label{sec-fp}

\subsection{Elementary properties}\label{sec-fp-basic}

\begin{proposition}\label{prop-classes}
\begin{enumerate}
\item[(1)] Every feebly compact space is a feebly $p$-space.
\item[(2)] Every feebly $p$-space is a feebly $q$-space.
\item[(3)] If every point of a space $X$ has a neighborhood which is a feebly compact subspace (in particular, if $X$ is locally compact),
then $X$ is a feebly $p$-space.
\end{enumerate}
\end{proposition}

\begin{proof}
(1) Put $\gamma_n=\{X\}$ for $n\in\om$ and use Proposition~\ref{prop-fc}.

(2) Let $(\gamma_n)$ be as in Definition~\ref{def-fp}, let $x\in X$ and choose $U_n\in\gamma_n$ with $x\in U_n$. The sets
$Q_n=U_0\cap\dots\cap U_n$ are open neighborhoods of $x$, $Q_{n+1}\subset Q_n$ and $Q_n\subset U_n$, so $(Q_n)$ is feebly compact by the
definition of a feebly $p$-space.

(3) Let $\gamma$ be the family of all open sets $U\subset X$ which are contained in a feebly compact subspace of $X$; it covers $X$. Put
$\gamma_n=\gamma$ for $n\in\om$. Let $x\in U_n\in\gamma$, let $V_n\subset U_n$ be nonempty open sets with $V_{n+1}\subset V_n$, and let
$M_n\subset V_n$ be nonempty open sets. Choose a feebly compact subspace $N\supset U_0$. The sets $T_n=\bigcup_{k\ge n}M_k\subset V_0\subset N$
are nonempty, open in $N$ and decreasing, so by Lemma~\ref{l-fc} applied to $N$ there is a point $p\in\bigcap_n\cl{T_n}=\Ls_nM_n$. Hence
$(V_n)$ is feebly compact.
\end{proof}

The class of regular feebly $p$-spaces is preserved by passing to open subspaces and to dense $G_\delta$-subspaces.

\begin{proposition}\label{prop-hered}
Let $X$ be a regular feebly $p$-space. Then every open subspace and every dense $G_\delta$-subspace of $X$ is a feebly $p$-space.
\end{proposition}

\begin{proof}
Let $(\gamma_n)_{n\in\om}$ be as in Definition~\ref{def-fp}.

\emph{Open subspaces.} Let $U\subset X$ be open, and let $\mathcal C$ be the family of open sets $C\subset X$ with $\cl{C}\subset U$; by
regularity, it covers $U$. The families $\gamma_n'=\{W\cap C:W\in\gamma_n,\ C\in\mathcal C\}$ are open covers of $U$. Let $x\in U$,
$x\in W_n\cap C_n\in\gamma_n'$, let $V_n\subset W_n\cap C_n$ be nonempty open sets with $V_{n+1}\subset V_n$, and let $M_n\subset V_n$ be
nonempty open sets. The sequence $(V_n)$ is feebly compact in $X$, since $x\in W_n\in\gamma_n$; hence $(M_n)$ accumulates to some
$p\in X$. By \eqref{eq-ls}, $p\in\cl{V_0}\subset\cl{C_0}\subset U$, so $(M_n)$ accumulates to $p$ in $U$.

\emph{Dense $G_\delta$-subspaces.} Let $H=\bigcap_n{E}_n$ be dense in $X$, where the sets ${E}_n$ are open. Let $\gamma_n^*$ be the family
of all open sets $W$ with $\cl{W}\subset W'$ for some $W'\in\gamma_n$, and let $\mathcal O_n$ be the family of open sets $O$ with $\cl{O}\subset {E}_n$; by
regularity, $\gamma_n^*$ covers $X$ and $\mathcal O_n$ covers $H$. The families $\gamma_n'=\{W\cap O\cap H:W\in\gamma_n^*,\ O\in\mathcal
O_n\}$ are open covers of $H$.
Let $x\in H$, $x\in W_n\cap O_n\cap H\in\gamma_n'$, let $V_n\subset W_n\cap O_n\cap H$ be nonempty open subsets of $H$ with
$V_{n+1}\subset V_n$, and let $M_n\subset V_n$ be nonempty open subsets of $H$. All closures below are taken in $X$. Choose open sets
$M_n'\subset X$ with $M_n=M_n'\cap H$ and put $F_n=\Int\cl{V_n}$. Since $H$ is dense, $M_n'\subset\cl{M_n}\subset\cl{V_n}$, hence
$\emptyset\ne M_n'\subset F_n$. The sets $F_n$ are open and decreasing. Since $W_n\in\gamma_n^*$, there is $W_n'\in\gamma_n$ with
$\cl{W_n}\subset W_n'$; then $x\in W_n'$ and $F_n\subset\cl{W_n}\subset W_n'$. By the definition of a feebly $p$-space applied to the sets
$W_n'$, the sequence $(F_n)$ is feebly compact, and $(M_n')$ accumulates to a point $p\in X$. By \eqref{eq-ls},
$p\in\cl{V_n}\subset\cl{O_n}\subset {E}_n$ for every $n$, so $p\in H$. Finally, if $W$ is a neighborhood of $p$ and $W\cap M_k'\ne\emptyset$,
then the nonempty open set $W\cap M_k'$ meets the dense set $H$, so $W\cap M_k\ne\emptyset$. Thus $(M_n)$ accumulates to $p$ in $H$.
\end{proof}

\subsection{\texorpdfstring{Feebly $p$-spaces versus $p$-spaces}{Feebly p-spaces versus p-spaces}}\label{sec-fp-vs-p}

The original definition of $p$-spaces by Arhangel'skii~\cite{Arhangelskii1963} uses the position of $X$ in
its Stone--\v{C}ech compactification, so it makes sense only for Tychonoff spaces.
Gruenhage~\cite[p.~441]{Gruenhage1984} defines $p$-spaces only for completely regular spaces and remarks that
the internal characterization of Burke \cite[Theorem~3.21]{Gruenhage1984} can be used as a definition.
Accordingly, for arbitrary spaces we call $X$ a \term{$p$-space} if there is a sequence $(\gamma_n)_{n\in\om}$ of open covers of $X$ such that
if $x\in X$ and $x\in U_n\in\gamma_n$ for $n\in\om$, then the set $K=\bigcap_n\cl{U_n}$ is compact and the sequence
$\bigl(\bigcap_{i\le n}\cl{U_i}\bigr)_{n\in\om}$ is an \term{outer network} for $K$, that is,
for every open $O\supset K$ there is $n$ with $\bigcap_{i\le n}\cl{U_i}\subset O$.
For Tychonoff spaces this is equivalent to the original definition \cite[Theorem~3.21]{Gruenhage1984}.
Recall also that a space $X$ is \term{pointwise countably complete} if there is a sequence $(\gamma_n)$ of open covers of $X$
such that $\bigcap_n\cl{A_n}\ne\emptyset$ for every decreasing sequence of nonempty sets $A_n\subset U_n\in\gamma_n$ with
$\bigcap_nU_n\ne\emptyset$ \cite[Section~2(3)]{Bouziad1996}; every $p$-space is pointwise countably complete (Proposition~\ref{prop-p-fp}; see also \cite{Bouziad1996}).

\begin{proposition}\label{prop-p-fp}
Every $p$-space is pointwise countably complete, and every pointwise countably complete space is a feebly $p$-space.
\end{proposition}
\begin{proof}
Let $(\gamma_n)$ be a sequence of open covers of a $p$-space $X$ as in the definition.
Let $A_n\subset U_n\in\gamma_n$ be nonempty decreasing sets with $x\in\bigcap_nU_n$.
Then $A_n\subset A_i\subset U_i$ for $i\le n$, so $\cl{A_n}\subset N_n:=\bigcap_{i\le n}\cl{U_i}$.
Pick $a_n\in A_n$. If $(a_n)$ had no cluster point in the compact set $K=\bigcap_n\cl{U_n}$,
every point of $K$ would have an open neighborhood containing $a_n$ for only finitely many $n$;
finitely many of them cover $K$, and their union $O$ is an open set containing $a_n$ for only finitely many $n$.
But $N_n\subset O$ for large $n$, a contradiction. A cluster point $a$ of $(a_n)$ lies in every $\cl{A_n}$, because $a_k\in A_n$ for $k\ge n$.
Hence $\bigcap_n\cl{A_n}\ne\emptyset$, and $X$ is pointwise countably complete.

Now let $X$ be pointwise countably complete, with covers $(\gamma_n)$. We show that the same covers witness that $X$ is a feebly $p$-space.
Let $x\in U_n\in\gamma_n$, let $V_n\subset U_n$ be decreasing nonempty open sets and let $M_n\subset V_n$ be nonempty open sets.
The sets $T_k=\bigcup_{n\ge k}M_n\subset U_k$ are nonempty and decreasing and $x\in\bigcap_kU_k$, so
$\Ls_nM_n=\bigcap_k\cl{T_k}\ne\emptyset$. Hence $(V_n)$ is feebly compact.
\end{proof}

The relations between the classes of spaces considered in this paper are summarized by the following diagram
(see Propositions~\ref{prop-classes} and~\ref{prop-p-fp}); we also note that a product of a $p$-space with a feebly $p$-space is a feebly $p$-space
(Proposition~\ref{prop-prod}).
\[
\begin{tikzcd}[column sep=small, row sep=small]
p\text{-spaces} \arrow[r, Rightarrow] & \parbox{3.6cm}{\centering pointwise countably complete} \arrow[r, Rightarrow] & \parbox{1.9cm}{\centering feebly $p$-spaces} \arrow[r, Rightarrow] & \parbox{1.9cm}{\centering feebly $q$-spaces}\\
 & & \parbox{1.9cm}{\centering feebly compact} \arrow[u, Rightarrow] &
\end{tikzcd}
\]

\begin{proposition}\label{prop-prod}
Let $X$ be a $p$-space and $Y$ a feebly $p$-space. Then $X\times Y$ is a feebly $p$-space.
\end{proposition}
\begin{proof}
Let $(\alpha_n)$ and $(\gamma_n)$ be sequences of open covers of $X$ and $Y$ as in the definitions of a $p$-space and of
a feebly $p$-space, respectively. We show that the covers $\delta_n=\{A\times B: A\in\alpha_n,\ B\in\gamma_n\}$ of $X\times Y$
witness that $X\times Y$ is a feebly $p$-space.

Let $(x,y)\in A_n\times B_n\in\delta_n$ for $n\in\om$, let $V_n\subset A_n\times B_n$ be nonempty open sets with
$V_{n+1}\subset V_n$, and let $M_n\subset V_n$ be nonempty open sets. Shrinking $M_n$ does not destroy accumulation points,
so we may assume that $M_n=P_n\times Q_n$ with $P_n\subset X$, $Q_n\subset Y$ nonempty and open.

Put $W_n=\mathrm{pr}_X(V_n)$ and $S_n=\mathrm{pr}_Y(V_n)$. These are nonempty open sets (projections are open),
and $P_n\subset W_n$, $Q_n\subset S_n$. The sets $S_n$ decrease and $S_n\subset B_n$, so $(S_n)$ is feebly compact by the choice of $(\gamma_n)$.
Hence $(Q_n)$ has an accumulation point in $Y$.
Moreover, $W_n\subset W_i\subset A_i$ for $i\le n$, so $\cl{W_n}\subset N_n:=\bigcap_{i\le n}\cl{A_i}$.
Let $K=\bigcap_n\cl{A_n}$; it is compact, and $(N_n)$ is an outer network for $K$.

Suppose that $(M_n)$ has no accumulation point in $X\times Y$. Fix $q\in Y$.
Every point of the compact set $K\times\{q\}$ has an open neighborhood $U\times V$ such that
$M_n\cap(U\times V)\ne\emptyset$ for only finitely many $n$. Choose finitely many of them covering $K\times\{q\}$ and let
$U^q$ be the union of their first coordinates and $V^q$ the intersection of their second coordinates.
Then $U^q\supset K$ is open, $V^q$ is a neighborhood of $q$, and $M_n\cap(U^q\times V^q)\ne\emptyset$ for only finitely many $n$.
Since $(N_n)$ is an outer network for $K$ and $P_n\subset\cl{W_n}\subset N_n$, there is $m$ with $P_n\subset U^q$ for $n\ge m$.
For such $n$ we get $M_n\cap(U^q\times V^q)=P_n\times(Q_n\cap V^q)$, which is nonempty whenever $Q_n\cap V^q\ne\emptyset$.
Hence $Q_n\cap V^q\ne\emptyset$ for only finitely many $n$, so $q$ is not an accumulation point of $(Q_n)$.
As $q\in Y$ was arbitrary, $(Q_n)$ has no accumulation point in $Y$, a contradiction.
Thus $(M_n)$ accumulates at a point of $X\times Y$, and $(V_n)$ is feebly compact.
\end{proof}

\begin{corollary}\label{cor-prod}
The product of a $p$-space and a feebly compact space is a feebly $p$-space.
\end{corollary}
\begin{proof}
Every feebly compact space is a feebly $p$-space by Proposition~\ref{prop-classes}(1); apply Proposition~\ref{prop-prod}.
\end{proof}

\begin{remark}
Regularity of the factors is not used in the proofs above; it is needed only if one wants the statement for
$p$-spaces in the original sense of Arhangel'skii.
\end{remark}

\subsection{\texorpdfstring{$\Gamma^{BM}_o$-nonmeager spaces}{Gamma-nonmeager spaces}}\label{sec-delta}\label{sec-gamma}
We use three classes of spaces defined by topological games in \cite{R2023}. Recall the Banach--Mazur game on a space $X$
(Section~\ref{sec-bm}): $U_0=X$, and at step $n\ge1$ player $\beta$ chooses a nonempty open set $V_n\subset U_{n-1}$, and then player
$\alpha$ chooses a nonempty open set $U_n\subset V_n$. The \term{modified Banach--Mazur game} on $X$ has the same rules, except that at step~$0$
player $\alpha$ chooses an arbitrary nonempty open set $U_0\subset X$. In both games a play determines the decreasing sequence $(V_n)_{n\ge1}$
of nonempty open sets, and $\bigcap_nV_n=\bigcap_nU_n$.

\begin{definition}\label{def-O}
Let $X$ be a space. A decreasing sequence $(V_n)_{n\ge1}$ of nonempty open subsets of $X$ satisfies \term{condition~\textup{(O)}} if
\begin{enumerate}
\item[(O)] for every sequence $(M_n)_{n\ge1}$ of nonempty open sets with $M_n\subset V_n$, the sequence $(M_n)$ has an accumulation point
in $\bigcap_{n\ge1}V_n$, that is, $\Ls_nM_n\cap\bigcap_nV_n\ne\emptyset$.
\end{enumerate}
\end{definition}

In the notation of \cite{R2023}, the sequences satisfying (O) form the family $\mathcal V_o(X)$. We note three simple properties of
condition~(O).
\begin{enumerate}
\item[(O1)] If $(V_n)$ satisfies (O), then $\bigcap_nV_n\ne\emptyset$.
\item[(O2)] Condition (O) depends only on the tail of the sequence: for every $k\ge1$, $(V_n)_{n\ge1}$ satisfies (O) if and only if
$(V_n)_{n\ge k}$ does.
\item[(O3)] If $U\subset X$ is open and $V_n\subset U$ for all $n$, then $(V_n)$ satisfies (O) in $U$ if and only if it satisfies (O) in $X$.
\end{enumerate}
Indeed, (O1) follows by taking $M_n=V_n$; (O2) holds because $\Ls_nM_n$ does not depend on finitely many terms and
$\bigcap_{n\ge1}V_n=\bigcap_{n\ge k}V_n$; and (O3) holds because the open subsets of $U$ are open in $X$ and every point of $\bigcap_nV_n$
lies in the open set $U$.

\begin{definition}[{\cite{R2023}}]\label{def-gamma}
A space $X$ is
\begin{enumerate}
\item[(1)] a \term{$\Gamma^{BM}_o$-space} if player $\alpha$ has a strategy in the Banach--Mazur game on $X$ such that every play which
follows it satisfies $\bigcap_nV_n=\emptyset$ or (O);
\item[(2)] \term{$\Gamma^{BM}_o$-Baire} if for every strategy of player $\beta$ in the Banach--Mazur game on $X$ there is a play which
follows it and satisfies (O);
\item[(3)] \term{$\Gamma^{BM}_o$-nonmeager} if for every strategy of player $\beta$ in the modified Banach--Mazur game on $X$ there is a
play which follows it and satisfies (O).
\end{enumerate}
\end{definition}

By Theorem~\ref{th-oxtoby} and (O1), $\Gamma^{BM}_o$-Baire spaces are Baire spaces, and $\Gamma^{BM}_o$-nonmeager spaces are nonmeager
\cite[Proposition~14]{R2023}. Every $\Gamma^{BM}_o$-Baire space is $\Gamma^{BM}_o$-nonmeager: a strategy of $\beta$ in the modified game,
applied to the first move $U_0=X$ of $\alpha$, is a strategy of $\beta$ in the Banach--Mazur game. Conversely, a Baire (nonmeager)
$\Gamma^{BM}_o$-space is $\Gamma^{BM}_o$-Baire ($\Gamma^{BM}_o$-nonmeager) \cite[Propositions~7 and~11]{R2023}. The relation between
$\Gamma^{BM}_o$-nonmeager and $\Gamma^{BM}_o$-Baire spaces is the same as between nonmeager and Baire spaces.

\begin{lemma}\label{l-gamma-local}
\begin{enumerate}
\item[(1)] Every open subspace of a $\Gamma^{BM}_o$-Baire space is $\Gamma^{BM}_o$-Baire.
\item[(2)] If every point of $X$ has an open neighborhood which is $\Gamma^{BM}_o$-Baire, then $X$ is $\Gamma^{BM}_o$-Baire.
\item[(3)] Every $\Gamma^{BM}_o$-nonmeager space contains a nonempty open $\Gamma^{BM}_o$-Baire subspace.
\item[(4)] Every homogeneous $\Gamma^{BM}_o$-nonmeager space, in particular every $\Gamma^{BM}_o$-nonmeager semitopological group, is
$\Gamma^{BM}_o$-Baire.
\end{enumerate}
\end{lemma}

\begin{proof}
(1) Let $O$ be an open subspace of a $\Gamma^{BM}_o$-Baire space $Y$, and let $t$ be a strategy of $\beta$ in the Banach--Mazur game on $O$.
It is also a strategy of $\beta$ in the Banach--Mazur game on $Y$, so some play which follows it satisfies (O) in $Y$, and hence in $O$
by (O3).

(2) Let $t$ be a strategy of $\beta$ in the Banach--Mazur game on $X$, and let $V_1$ be its first move. Take $x\in V_1$ and an open
$\Gamma^{BM}_o$-Baire neighborhood $U$ of $x$; by (1), $O=V_1\cap U$ is $\Gamma^{BM}_o$-Baire. Consider the strategy of $\beta$ in the
Banach--Mazur game on $O$ whose first move is $O$ and which then follows $t$ (the moves of $\alpha$ in $O$ are admissible answers to $V_1$ in
the game on $X$). Some play of this strategy satisfies (O) in $O$, hence in $X$ by (O3); the corresponding play following $t$ differs from it only
in the first move of $\beta$, hence it also satisfies (O) by (O2).

(3) Suppose that no nonempty open subspace of $X$ is $\Gamma^{BM}_o$-Baire. Then for every nonempty open $U\subset X$ there is a strategy
$t_U$ of $\beta$ in the Banach--Mazur game on $U$ such that no play following it satisfies (O). In the modified Banach--Mazur game on $X$,
let $\beta$ answer the first move $U_0$ of $\alpha$ by following $t_{U_0}$. By (O3), no play following this strategy satisfies (O), so $X$ is
not $\Gamma^{BM}_o$-nonmeager.

(4) By (3), some nonempty open set $U$ is $\Gamma^{BM}_o$-Baire. By homogeneity every point has an open neighborhood homeomorphic to $U$,
and (2) applies.
\end{proof}

\begin{remark}\label{rem-ACK}
$\Gamma^{BM}_o$-Baire spaces were in fact considered in \cite{ACK2010} and \cite{R2022ci}. In \cite{ACK2010}, a decreasing sequence
$(H_n)$ of nonempty open sets has property $\Pi$ if $\bigcap_n\cl{W_n}\ne\emptyset$ for every decreasing sequence of nonempty open sets
$W_n\subset H_n$, and player $\alpha$ wins a play of the Banach--Mazur game in the game $G_\Pi$ if $(V_n)$ has property $\Pi$. For a regular
space $X$, the following conditions are equivalent: (1)~$X$ is $\Gamma^{BM}_o$-Baire; (2)~$\beta$ has no winning strategy in $G_\Pi$.
Indeed, (1)$\Rightarrow$(2) is obvious, since for decreasing $W_n\subset V_n$ every accumulation point of $(W_n)$ lies in $\bigcap_n\cl{W_n}$.
Conversely, let $t$ be a strategy of $\beta$ in the Banach--Mazur game. Using regularity, let $\beta$ answer, instead of a move $V_n$
prescribed by $t$, a nonempty open set $V_n'$ with $\cl{V_n'}\subset V_n$, and let $t$ be applied to the actual moves $U_n\subset V_n'$ of
$\alpha$. By (2), there is a play of this new strategy in which $(V_n')$ has property $\Pi$; the sets $(V_n,U_n)$ form a play which follows
$t$. If $M_n\subset V_n$ are nonempty open sets, then $T_k=\bigcup_{n>k}M_n\subset V_{k+1}\subset V_k'$ are decreasing, so
$\Ls_nM_n=\bigcap_k\cl{T_k}\ne\emptyset$, and $\Ls_nM_n\subset\bigcap_k\cl{V_{k+1}'}\subset\bigcap_kV_k$ by
\eqref{eq-ls}. Thus the play $(V_n,U_n)$ follows $t$ and satisfies (O), and $X$ is $\Gamma^{BM}_o$-Baire.

In \cite{R2022ci}, player $\alpha$ wins a play of the Banach--Mazur game in the game $G(BM,b)$ if the sequence $(U_n)$ of his moves is
\emph{bounded}, that is, for every locally finite family $\gamma$ of open sets some $U_n$ meets only finitely many members of $\gamma$. For
decreasing sequences of nonempty open sets this means exactly that every sequence of nonempty open sets $M_n\subset U_n$ has an
accumulation point (if a sequence $(M_n)$ has no accumulation point, then no set occurs in it infinitely many times, and the family
$\{M_n:n\ge1\}$ is locally finite). Since $V_{n+1}\subset U_n\subset V_n$, the same argument shows that for a regular space $X$ the
absence of a winning strategy of $\beta$ in $G(BM,b)$ is also equivalent to $X$ being $\Gamma^{BM}_o$-Baire.
\end{remark}

\begin{proposition}\label{prop-gamma-space}
Every regular feebly $p$-space is a $\Gamma^{BM}_o$-space.
\end{proposition}

\begin{proof}
Let $(\gamma_n)_{n\in\om}$ be as in Definition~\ref{def-fp}. We define a strategy for $\alpha$. Suppose that $\beta$ has chosen a nonempty
open set $V_n$, $n\ge1$. Take a point $x\in V_n$ and $W\in\gamma_{n-1}$ with $x\in W$. By regularity, there is an open set $U_n\ni x$ with
$\cl{U_n}\subset V_n\cap W$; player $\alpha$ chooses this $U_n$. Thus
\begin{equation}\label{eq-alpha}
\cl{U_n}\subset V_n\quad\text{and}\quad U_n\subset W_{n-1}\text{ for some }W_{n-1}\in\gamma_{n-1}.
\end{equation}
Consider a play which follows this strategy. Since $V_{n+1}\subset U_n$, we have $\cl{V_{n+1}}\subset V_n$ for all $n\ge1$. Suppose that
$\bigcap_nV_n\ne\emptyset$ and take $x\in\bigcap_nV_n$. By \eqref{eq-alpha}, $x\in V_{m+2}\subset U_{m+1}\subset W_m\in\gamma_m$ for
$m\in\om$. The sets $V'_m=V_{m+2}$ are nonempty, open and decreasing, and $V'_m\subset W_m$; hence the sequence $(V'_m)_{m\in\om}$ is
feebly compact. Let $M_n\subset V_n$ be nonempty open sets. Since $M_{m+2}\subset V'_m$, the sequence $(M_n)$ accumulates to a point $p$.
For every $k\ge1$ we have $M_n\subset V_{k+1}$ for $n\ge k+1$, so $p\in\cl{V_{k+1}}\subset V_k$ by \eqref{eq-ls}. Thus $(M_n)$ has the
accumulation point $p\in\bigcap_kV_k$, so $(V_n)$ satisfies (O), and the strategy is as required.
\end{proof}

\begin{proof}[Proof of Proposition~\ref{prop-fp-delta}]
Let $X$ be a regular Baire feebly $p$-space. By Proposition~\ref{prop-gamma-space}, $X$ is a $\Gamma^{BM}_o$-space. A nonmeager
$\Gamma^{BM}_o$-space is $\Gamma^{BM}_o$-nonmeager \cite[Propositions~7(1) and~11]{R2023}, and every $\Gamma^{BM}_o$-nonmeager space is
$\Gamma^{OD}_{o,l}$-nonmeager \cite[Proposition~9(1)]{R2023}. (In the notation of \cite{R2023}, the sequences satisfying (O) form the family
$\mathcal V_w(X,\mathcal N_o,\mathcal W_l)$.) Finally, $\Gamma^{OD}_{o,l}$-nonmeager spaces are
$\Delta$-nonmeager \cite[Propositions~5.1 and~6.1(1)]{R2025}.
\end{proof}

The proof of Proposition~\ref{prop-fp-delta} shows that every $\Gamma^{BM}_o$-nonmeager space is $\Delta$-nonmeager. Another useful
property of $\Gamma^{BM}_o$-nonmeager spaces is the existence of feebly $q$-points.

\begin{lemma}\label{l-gamma-q}
Every $\Gamma^{BM}_o$-nonmeager space contains a feebly $q$-point. In a $\Gamma^{BM}_o$-Baire space the set of feebly $q$-points is dense.
\end{lemma}

This is essentially \cite[Proposition~12]{R2023}: the points of $\bigcap_nV_n$ for sequences $(V_n)$ satisfying (O), called
$q_o$-points in \cite{R2023}, are feebly $q$-points. For completeness we give the short proof.

\begin{proof}
Let $X$ be $\Gamma^{BM}_o$-nonmeager, and let $\beta$ always answer with the set just chosen by $\alpha$. There is a play which follows this
strategy and satisfies (O). By (O1), there is a point $x\in\bigcap_nV_n$. Then $(V_n)$ is a sequence of open neighborhoods of $x$, and by
(O) it is feebly compact; hence $x$ is a feebly $q$-point. If $X$ is $\Gamma^{BM}_o$-Baire and $W\subset X$ is a nonempty open set,
apply the same argument to the strategy of $\beta$ in the Banach--Mazur game whose first move is $V_1=W$; then $x\in W$.
\end{proof}

\subsection{Groups}

\begin{proof}[Proof of Theorem~\ref{th-main-gamma}]
Let $G$ be a regular separable $\Gamma^{BM}_o$-nonmeager semitopological group. By Lemma~\ref{l-gamma-local}(4), $G$ is
$\Gamma^{BM}_o$-Baire, and hence a Baire space. By Lemma~\ref{l-gamma-q}, $G$ has a feebly $q$-point; since translations are
homeomorphisms, every point of $G$ is a feebly $q$-point. By Theorem~\ref{th-main-scf} with $X=Y=Z=G$, the multiplication on $G$ is quasicontinuous. As shown in
the proof of Proposition~\ref{prop-fp-delta}, $G$ is $\Delta$-nonmeager, so $G$ is a topological group by Theorem~\ref{th-delta-stg}.
\end{proof}

\begin{proof}[Proof of Theorem~\ref{th-main-pspace}]
Let $G$ be a regular separable Baire semitopological group which is a feebly $p$-space. By Proposition~\ref{prop-gamma-space}, $G$ is a
$\Gamma^{BM}_o$-space, and a nonmeager $\Gamma^{BM}_o$-space is $\Gamma^{BM}_o$-nonmeager \cite[Propositions~7(1) and~11]{R2023}. Apply
Theorem~\ref{th-main-gamma}.
\end{proof}

Theorem~\ref{th-main-group} is also a special case of Theorem~\ref{th-main-gamma}: a regular feebly compact space is a Baire space
(Lemma~\ref{l-baire}) and a feebly $p$-space (Proposition~\ref{prop-classes}(1)).

\subsection{Local versions}

\begin{remark}[Local versions]\label{rem-local}
A space is \term{locally separable} if each of its points has a separable open neighborhood, and \term{locally feebly $p$} if each of its
points has an open neighborhood which is a feebly $p$-space. The main results of the paper remain valid in the following local forms:
\begin{itemize}
\item in Theorems~\ref{th-main-group}, \ref{th-main-scf}, \ref{th-main-scf+}, \ref{th-main-pspace} and~\ref{th-main-gamma},
``separable'' can be replaced by ``locally separable'';
\item in Proposition~\ref{prop-fp-delta} and Theorem~\ref{th-main-pspace}, ``feebly $p$-space'' can be replaced by ``locally feebly
$p$-space''.
\end{itemize}
Indeed:
\begin{enumerate}
\item[(i)] \emph{Local separability.} In the proof of Theorem~\ref{th-main-scf+}, separability of $X$ is used only in Step~1, to find a
countable dense subset of the open neighborhood $U_*$ of $x_*$. If $S$ is a separable open neighborhood of $x_*$, we may choose
$U_*\subset S\cap U_*'$; then $U_*$ is separable as an open subspace of the separable space $S$, and the rest of the proof is unchanged.
Hence Theorems~\ref{th-main-scf+} and~\ref{th-main-scf} hold for locally separable Baire spaces $X$, and the proofs of
Theorems~\ref{th-main-group}, \ref{th-main-pspace} and~\ref{th-main-gamma} in Sections~\ref{sec-group} and~\ref{sec-fp} apply
without changes. Since translations of a semitopological group are homeomorphisms, a semitopological group is locally separable if and only if it has a nonempty
separable open subset. In fact, in Theorems~\ref{th-main-scf+} and~\ref{th-main-scf} it suffices that the set of points of $X$ having a
separable open neighborhood is dense in $X$. Indeed, after Step~1 the point $x_*$ is not used, and the proof only needs a nonempty
separable open set $U_*\subset U_*'$ with $f(U_*\times\{y_*\})\subset A$; one can take for $U_*$ any nonempty separable open subset of an open neighborhood
$U$ of $x_*$ with $U\subset U_*'$ and $f(U\times\{y_*\})\subset A$. For semitopological groups this condition is equivalent to local
separability.
\item[(ii)] \emph{Every locally feebly $p$-space $X$ is a feebly $q$-space.} Let $x\in X$ and let $O$ be an open neighborhood of $x$ which is a
feebly $p$-space. By Proposition~\ref{prop-classes}(2), $x$ has a sequence $(Q_n)$ of open neighborhoods which is feebly compact in $O$. If
$M_n\subset Q_n$ are nonempty open subsets of $X$, then they are open in $O$, so $(M_n)$ has an accumulation point in $O$, which is also an
accumulation point in $X$. Hence $(Q_n)$ is feebly compact in $X$, and Theorem~\ref{th-main-scf} applies with this space in place of $Y$.
\item[(iii)] \emph{Every regular locally feebly $p$-space is a $\Gamma^{BM}_o$-space.} Every point has an open neighborhood which is a
regular feebly $p$-space, hence a $\Gamma^{BM}_o$-space by Proposition~\ref{prop-gamma-space}, and a space in which every point has a
$\Gamma^{BM}_o$-neighborhood is a $\Gamma^{BM}_o$-space \cite[Theorem~3(5)]{R2023}. The proof of Proposition~\ref{prop-fp-delta} then
applies without changes, and so does the proof of Theorem~\ref{th-main-pspace}, by (ii).
\end{enumerate}
Similarly, a semitopological group is locally feebly $p$ if and only if its identity has an open neighborhood which is a feebly
$p$-space.
\end{remark}

A space is \term{locally feebly compact} if each of its points has a neighborhood which is a feebly compact subspace.

\begin{corollary}\label{cor-local}
Every regular locally separable locally feebly compact semitopological group is a topological group.
\end{corollary}

\begin{proof}
Let $G$ be such a group. By Proposition~\ref{prop-classes}(3), $G$ is a feebly $p$-space. Let us show that $G$ is a Baire space. Every
point $x\in G$ has a neighborhood $N$ which is a feebly compact subspace; $N$ is regular, hence a Baire space by Lemma~\ref{l-baire}. The
set $\Int N$ is an open subspace of $N$, hence also a Baire space. Thus every point of $G$ has an open neighborhood which is a Baire space,
and therefore $G$ is a Baire space. By Remark~\ref{rem-local}(i), Theorem~\ref{th-main-pspace} holds for locally separable groups, so $G$
is a topological group.
\end{proof}

\section{Open problems}\label{sec-problems}

Cases (b) and (c) of Question~\ref{q-T} remain open. Since every Fr\'echet--Urysohn space has countable tightness and is a $k$-space,
the Fr\'echet--Urysohn case is the weakest special case of both; we include it explicitly.

\begin{question}
Let $G$ be a regular feebly compact semitopological group which is a Fr\'echet--Urysohn space, a space of countable tightness, or
a $k$-space. Is $G$ a topological group?
\end{question}

The corresponding question on separately continuous maps is the following.

\begin{question}\label{q-qc}
Let $X$ and $Y$ be regular feebly compact spaces, where $X$ is a Fr\'echet--Urysohn space, a space of countable tightness, or a
$k$-space, and let $Z$ be a regular space. Is every separately continuous map $f\colon X\times Y\to Z$ strongly quasicontinuous?
\end{question}

For Tychonoff groups $G$ and Tychonoff spaces $X$, $Y$, $Z$ both questions have positive answers \cite{R1994,R2022}.

\end{document}